\documentclass[11pt,reqno]{amsart}
\usepackage{amssymb,latexsym}
\usepackage{graphicx}
\usepackage{fancyhdr,amssymb}
\usepackage{url}		

\theoremstyle{plain}
\numberwithin{equation}{section}

\begin{document}
\fancyhead{}
\renewcommand{\headrulewidth}{0pt}
\fancyfoot{}
\fancyfoot[LE,RO]{\medskip \thepage}

\setcounter{page}{1}

\title[Fibonacci Numbers and Identities ]{Fibonacci Numbers and Identities}
\author{Cheng Lien Lang}
\address{Department Applied of Mathematics\\
                I-Shou University\\
                Kaohsiung, Taiwan\\
                Republic of China}
\email{cllang@isu.edu.tw}
\thanks{}
\author{Mong Lung Lang}
\address{     }
\email{lang2to46@gmail.com}

\begin{abstract}
By investigating  a recurrence relation about functions, we first give
alternative proofs  for various identities on Fibonacci numbers and Lucas numbers, and then,
 make certain well known identities {\em visible} via certain trivalent graph
  associated to the recurrence relation.

\end{abstract}

\maketitle

\vspace{-1cm}

\section{Introduction}
A function $x(n)$ defined over $\Bbb N \cup \{0\}$  is called an $\mathcal F$-function
 if $x(n)$ satisfies  the following recurrence relation.
 $$ x(n+3) = 2(x(n+2) + x(n+1)) - x(n).\eqno(1.1)$$
 One sees easily that the following are ${\mathcal F}$-functions :
   $$x(n) = (-1)^n, \,\,x(n) = F_{n }^2,\,\,
   x(n) = F_{n+r}F_{n},\,\,x(n) = F_{2n}, \eqno(1.2)$$
$$x(n) = L_{n}^2,\,\,
   x(n) = L_{n+r}L_{n},\,\,x(n) = L_{2n}, x(n) = F_nL_{n+r},\eqno(1.3)$$
    where $r \in \Bbb Z$ and $F_n,\, L_n$ are  the $n$-th Fibonacci and Lucas numbers
     respectively.
    Note that sum and difference
      of $\mathcal F$-functions are  $\mathcal F$-functions.
      A search of the literature turns out that a great deal of identities involve $\mathcal F$-functions
       only.
       For instance, all the  terms in the Cassini's identity $F_{n-1}F_{n+1}-F_n^2 = (-1)^n$ are
        $\mathcal F$-functions. For our convenience, we shall call such an identity $\mathcal F$-identity.
      Identities with other recurrence relations are less frequent.
      The main purpose of this article is to
      avoid the intricate case-by-case analysis, thereby obtaining
      a unified proof of the $\mathcal F$-identities.
      Since these identities involve $\mathcal F$-functions only, our proof will make usage of  (1.1)
       and (i), (ii) of the following only.
        \begin{enumerate}
    \item[(i)] (1.2) and (1.3) of the above are $\mathcal F$-functions,
    \item[(ii)]  $F_{n+2} = F_{n+1}+ F_n$, $F_{-m} = (-1)^{m+1}F_m,$
     $L_{n+2} = L_{n+1}+ L_n$, $L_{-m}= (-1)^mL_m$.
\end{enumerate}
Note that our proof can be applied easily to all $\mathcal F$-identities.
  Identities involve  other recurrence relations such as (iii) and (iv) of the
    following will be discussed in section 6.

   \begin{enumerate}
   \item[(iii)] $ A(n+2)= -A(n+1)+  A(n) $,
   \item[(iv)] $A(n+3) = -A(n+2)+A(n+1) +A(n)$.
   \end{enumerate}

The rest of the article is organised as follows. In section 2 we give some basic properties about
 $\mathcal F$-function. Section 3 gives alternative  proofs of the well known Catalan's identity
  and Melham's identity.
 Section 4 lists a few more identities (including d'Ocagne's, Tagiure's  and
 Gelin-Ces\`{a}ro
  identities) that can be proved by applying our technique presented in
  section 3. They are the $\mathcal F$-identities involved functions we listed in (1.2) and (1.3). In other
   words, they use functions in (1.2) and (1.3) as building blocks (see Lemma 2.1). Since product of $\mathcal
    F$-functions are not necessarily $\mathcal F$-functions, our idea cannot be applied to all
     identities (see Appendix B). Section 5 is devoted to the possible visualisation of identities via the
      recurrence relation (1.1). After all, there is nothing to prove if one cannot {\em  see} the
       identities in the first place. The last section gives a very brief discussion about
        identities that involve other recurrence relations.

    \section{Basic Properties about $\mathcal F$-functions}

\medskip
\noindent {\bf Lemma 2.1.} {\em Functions defined in $(1.2)$ and $(1.3)$ are $\mathcal F$-functions.
 Let $A(n)$ and $B(n)$ be $\mathcal F$-functions and let $r_0\in \Bbb Z$ be fixed. Then
  $X(n) = A(n +r_0), Y(n) = r_0 A(n)$ and $A(n)\pm B(n)$ are $\mathcal F$-functions. }
\medskip

\noindent {\em Proof.}  Let $x(n)$ be given as in (1.2) or (1.3). To show $x(n)$ is an $\mathcal F$-functions, it suffices
 to show that $x(n+3) = 2(x(n+2)+x(n+1))-x(n)$, which can be verified easily.   \qed

 \medskip
\noindent {\bf Lemma 2.2.} {\em
 Let $A(n)$ and $B(n)$ be $\mathcal F$-functions. Then
   $A(n)= B(n)$ if and only if
  $A(k)= B(k)$ for
$k = 0,1$ and $2$. }
\medskip

Note that the fact that the following functions are $\mathcal F$-functions has been treated
 as identities in the literature ([HB], identities (31)-(34) of [L]).

\medskip
\noindent {\bf Example 2.3.} By Lemma 2.1,
$F_{n+3}^2$, $F_{n+3}F_{n+4}$, $L_{n+3}^2$ and $ F_{n+3}L_{n+3}$
 are $\mathcal F$-functions.

  \medskip

 The following lemma is straightforward and will be used in sections 3 and 6.

   \medskip
   \noindent {\bf Lemma 2.4.} {\em Let $A(n)$ and $  B(n)$ be functions defined over $\Bbb N\cup \{0\}$.
    Suppose that both $A(n)$ and $B(n)$ satisfy either
    \begin{enumerate}
    \item[(i)] $x(n+2) = -x(n+1)+ x(n)$, or
    \item[(ii)] $x(n+3) = -2x(n+2) +2x(n+1) +x(n)$.
    \end{enumerate}
    Then $A(n) = B(n)$ if and only of $A(k)=B(k)$ for $k =0, 1,2$.}

\subsection {Discussion}  Let $\{x_n\}$ be a sequence  satisfies the recurrence relation $x_{n+2} = x_
 {n+1} +x_{n}$. Then $A(n) = x_{2n}$, $B(n) = x_nx_{n+r}, C(n) = x_n^2$ are $\mathcal F$-functions.

\section{ An alternative proof for Catalan's Identities}
In [H], Howard studied generalised Fibonacci sequence and  proved that Calatan's identity is equivalent to an identity discovered and proved by
 Melham ([M]) (see section 3.1). We give in the following our alternative proof which uses Lemmas 2.1 and
  2.2 only.

 \medskip

 \noindent {\bf Lemma 3.1.} $ 4(-1)^{3-r} + F_{r+3}F_{r-3} - F_{r}^2 =0.$

 \medskip
 \noindent {\em Proof.}
 We assume that $r \ge 0$. The case $r \le 0$ can be dealt with similarly. Let $A(r) = 4(-1)^{3-r} + F_{r+3}F_{r-3} - F_{r}^2 $. By Lemma 2.1, $A(r)$
  is an $\mathcal F$-function.
  By Lemma 2.2, $A(r) = 0$ for all $r$. This completes the proof of the Lemma. \qed
 \medskip

 \noindent {\bf Remark.} Any $\mathcal F$-identity with one variable $n$
 can be proved by applying the proof of
  Lemma 3.1. For instance, the Cassini's identity.

 \medskip
\noindent {\bf Theorem 3.2 (Catalan's Identity).} {\em
$F_n^2 -F_{n+r}F_{n-r} = (-1)^{n-r} F_r^2$.}

\medskip
\noindent {\em Proof.}   Recall first that
 $F_{-m} = (-1)^{m+1}F_m$. As a consequence, we may assume without loss of generality that
 $ n , r \ge 0$. As a matter of fact, we may assume that $n \ge 1$ as the case $n =0$ is trivial.
Let $A(n) =F_n^2 -F_{n+r}F_{n-r}$, $B(n)= (-1)^{n-r} F_{r}^2$.
 By Lemma 2.1,  both $A(n)$ and $B(n)$
  are $\mathcal F$-functions  (in $n$).
 Note that $$A(1) = 1-F_{1+r}F_{1-r}, \,\,A(2) = 1-F_{2+r}F_{2-r},\,\, A(3) = 4-F_{3+r}F_{3-r},
 \eqno (3.1)$$
  and $$B(1) = (-1)^{1-r}F_{r}^2, \,\,B(2) = (-1)^{2-r}F_{r}^2,\,\,
   B(3) = (-1)^{3-r}F_{r}^2.\eqno(3.2)$$

 \medskip
 \noindent  One sees easily that if $r \le 3$, then $A(n) = B(n)$ for $n=1,2,3$.
   By Lemma 2.2, we have $A(n) = B(n)$. Hence the theorem is proved.
    We shall therefore assume that $r\ge 4$.
  Recall that
 $F_{-m} = (-1)^{m+1}F_m$. This allows us to rewrite $A(3)$  into
 $A(3) = 4 +(-1)^{3-r} F_{r+3}F_{r-3}.$
 Hence $$A(3)-B(3) = 4 + (-1)^{3-r}F_{r+3}F_{r-3} - (-1)^{3-r}F_{r}^2.\eqno(3.3)$$

\medskip
\noindent By Lemma 3.1, $A(3) = B(3)$.
 One can show similarly that $A(1) = B(1)$ and $A(2) = B(2)$.
  Applying Lemma 2.2, we have
   $A(n) = B(n)$
for all $n $.  This completes the proof of the theorem. \qed

\medskip
\noindent {\bf Remark.} Any $\mathcal F$-identity with two variables $m$ and $n$
 can be proved by applying the proof of Theorem 3.2.

 \subsection{Melham's Identity.}  In [M], Melham proved among some very
  general results the identity  $F_{n+r+1}^2 +F_{n-r}^2
  = F_{2r+1}F_{2n+1}$. We shall give our alternative proof as follows.
  Denoted by $A(n)$ and $B(n)$ the left and right hand side of the identity.
   By Lemma 2.1,
   $A(n)$ and $B(n)$  are $\mathcal F$-functions (in $n$).
    The cases $n = 0,1$ and 2 of the above identity are given by
    $$F_{r+1}^2 + F_{-r}^2= F_{2r+1} F_1,\,\,
      F_{r+2}^2 + F_{1-r}^2= F_{2r+1} F_3,\,\,
      F_{r+3}^2 + F_{2-r}^2= F_{2r+1} F_5.\eqno(3.4)$$

      \noindent  By Lemma 2.1, the functions in (3.4) are $\mathcal F$-functions
       (in $r$)
      and that the identities  can be verified by applying Lemma 2.2. Consequently,
       we have $A(0)= B(0), A(1)=B(1)$ and $A(2)= B(2)$. By Lemma 2.2, we have
        $A(n) = B(n)$ for all $n$.

\subsection{Discussion.}  Our method can be generalised to functions such as $x(n)
 = F_n^3, y(n) = F_{3n}$ which satisfy the recurrence relation (see Appendix B)
  $$x(n+4) = 3x(n+3)+6x(n+2) -3x(n+1) -x(n).\eqno(3.5)$$
   \section{ More Identities}
 The purpose of this section is to list a few identities we found in the literature
  that can be proved by applying Lemmas 2.1 and 2.2
\subsection{ d'Ocagne's Identity}  The proof we presented in section 3 can be applied
 to all $\mathcal F$-identities.  A search
   of the literature turns out that there are many such identities. However, as
    the identities may be described in different manner, it is important to get
     equivalent forms of the  identities.  Take d'Ocagne's identity
  for instance. It is, in many occasion, given as follows (see [W]).

  $$F_mF_{n+1} - F_nF_{m+1} = (-1)^n F_{m-n}.\eqno(4.1)$$

\medskip
\noindent The very first look of the left and right hand side does not reveal the fact
 that they are $\mathcal F$-function. However, one has the following.
 Let $r= m-n$. Then (4.1) can be rewritten as
  $$F_{n+1}F_{n+r}-F_nF_{n+r+1} = (-1)^nF_{r},\eqno(4.2)$$

 \medskip
 \noindent
 where both the left and right hand side in (4.2) are $\mathcal F$-functions in terms of $n$.
 One may now apply the proof of Theorem 3.2 to give a proof of (4.2). As the proof is
  identically the same, we will not include it here.
\subsection{Some more identities}  A search of the literature turns up that there are many
   identities can be verified by Lemmas 2.1 and 2.2
   (for instance, out of the 44 identities given by Long [L], 35 of them
     involve $\mathcal F$-functions). We shall list a few which we pick mainly from [W]
    ($(c1)$-$(c8)$, $(c11)$, $(c12)$, $(d1)$-$(d6)$).

{\small
$\begin{array} {lllrrr}

\\
(c1) &F_{n+a}F_{n+b}-F_nF_{n+a+b}= (-1)^nF_aF_b    & :&   F_{2n}= F_{n+1}^2-F_{n-1}^2 & (d1)\\

\\

(c2) &  F_{n+1}^2 = 4F_n F_{n-1}+F_{n-2}^2    & : &  F_{2n+1}= F_{n+1}^2+F_{n}^2  & (d2) \\

\\

(c3) & L_n^2 -5F_n^2 = 4(-1)^n      & :&   F_{n+2}F_{n-1}= F_{n+1}^2 -F_{n}^2& (d3)\\

\\

(c4) & F_{n-1}F_{n+1}-F_n^2 = (-1)^{n-1}     & :&   F_{n}^2 - F_{n-2}F_{n+2} = (-1)^n& (d4)\\

\\

(c5) &F_mF_n = \frac{1}{5}(L_{m+n} -(-1)^nL_{m-n})    & :&   \sum_{n=1}^n F_n^2 = F_n F_{n+1}& (d5)\\

\\
(c6) & F_n^2 = \frac{1}{5}( L_{2n}-2(-1)^n) & :&   F_n^2 -F_{n-1}F_{n+1} = (-1)^{n-1}& (d6)\\

\\

(c7) & F_{n+m} = F_{n-1}F_m +F_nF_{m+1} & :&  F_n F_{n+3} = F_{n+1} F_{n+2} +(-1)^{n-1}& (d7)\\

\\

(c8)&F_{m+n} = \frac{1}{2}(F_mL_n+L_mF_n)   & :&  F_n^2 -F_{n-1}^2 = F_n F_{n-1} +(-1)^{n-1}& (d8)\\
\\
(c9) &F_{n+k+1}^2+F_{n-k}^2 =F_{2k+1}F_{2n+1}& :&   F_{2n+1}+(-1)^n =F_{n-1}F_{n+1} +F_{n+1}^2& (d9)\\
\\

(c10) &L_{n+k+1}^2+L_{n-k}^2 =5L_{2k+1}L_{2n+1}     &:& L_{2n+1}-F_{n+1}^2- A=(-1)^{n-1}& (d10)\\
\\
(c11) &F_n^2 +(-1)^{n+r-1} F_r^2 = F_{n-r}F_{n+r}   &:& L_{n-1}^2-F_{n-4}F_n-F_nF_{n+1}=F_{n-2}^2& (d11)\\

\\
(c12) & F_n^4 -F_{n-2}F_{n-1}F_{n+1}F_{n+2} = 1  &:&  F_{2n+1} = F_{n+3}F_n-F_{n+1}F_{n-1}& (d12)\\

\end{array}
$}

\bigskip
\noindent where  $A = (L_n^2-F_{n-3}F_{n+1})+F_{2n-2}$ and $L_n$ is the $n$-th Lucas number.
 (d10) and (d11) are less standard. We decide to include them in the table as they are
  {\em visible} via certain trivalent graph (see section 5).

\medskip

\noindent {\em Proof.} We note first that functions in $(c12)$ are not $\mathcal F$-functions but
 the identity can be proved by applying $(c4)$ and $(d4)$.
In $(c5), \,(c7)$ and $(c8)$, one needs to do a rewriting to see that
 the functions are $\mathcal F$-functions (let $m+n=r$).
   One may now apply Lemmas 2.1 and 2.2 and our proof presented in Theorem
    3.2  to verify these identities. \qed

\medskip
\noindent {\bf Remark.}
As identities may be described differently, the technique of rewriting identities into
 equivalent forms is crucial (see (4.1) and (4.2)).
$(c9)$ and $(c10)$
were first found and proved by Melhem ([M1]).
$(c12)$ is the Gelin-Ces\`{a}ro identity.

\subsection{Discussion.} Note that in our proof, we do not use any existing identities such as
 Binet's formula or any identities listed in [W] except (1.1), (1.2) and (1.3) of this
  article, which is what we promise in our introduction. Note also that one has to apply Lemma
  2.2 three times to prove identity $(c1)$, known as the Tagiuri's identity.

 \section {How Far can (1.1) go ?}

 We have demonstrated that the recurrence relation (1.1) can be used to verified various
  identities. In this section, we will present a trivalent graph (see the
   graph given in Appendix A) which is closely related
   to (1.1) that enables us to {\em visualise}  identities in the following.

 Let  $e_3,e_2,e_1$ be  vectors
 placed in the following
 trivalent graph and let   $e_4$ be the vector  given by
$ e_4 = 2(e_3+e_2)-e_1.$ Such a vector  $e_4$ is said to be
 $\mathcal F$-generated by $e_3$, $e_2$, and $e_1$ (in this order).
For our convenience, we use the following notation for the vector
 $e_4$ :

$$ e_4 = \left < e_3,e_2 :e_1\right > = 2(e_3+e_2)-e_1.\eqno(5.1)$$

\medskip
\noindent
Note that $ \left < e_3,e_2 :e_1\right >  \ne  \left < e_2,e_1 :e_3\right > \ne
 \left < e_1,e_3:e_2\right > $.
 Note also that (5.1) can be viewed as a generalisation (in the form of vectors) of
  the recurrence relation (1.1). We may construct an infinite sequence of vectors
   given as follows.
 $$e_1,e_2,e_3, e_4 = 2(e_3+e_2)-e_1, \cdots, e_{n+1} = 2(e_n +e_{n-1})-e_{n-1}, \cdots .    \eqno(5.2)$$

\bigskip
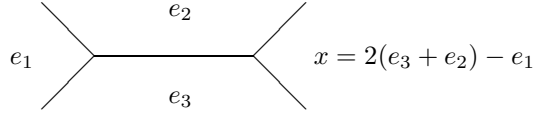
\begin{center}
 \begin{picture}(60,20)

 \multiput(2,0)(30,0){1}{\line(1,0){60}}

\multiput(2,0)(60,0){1}{\line(-1,1){20}}
 \multiput(2,0)(60,0){1}{\line(-1,-1){20}}

 \multiput(82,-20)(60,0){1}{\line(-1,1){20}}
 \multiput(82,20)(30,30){1}{\line(-1,-1){20}}

\put(-30, -3){\small  $e_1$}
\put(85, -3){\small $x = 2(e_3+e_2)-e_1$}
\put(30, 15){\small $ e_2$}
\put(30, -18){\small $e_3$}

 \end{picture}
\end{center}

\vspace{1cm}

\noindent
Denoted by $F(e_1,e_2,e_3)$ the above sequence.  In the
 case
  $\{e_1, e_2, e_3\}$ is  the canonical basis of $\Bbb R^3$,
the first nine vectors are given as follows.
 $e_1= (1,0,0),e_2=(0,1,0), e_3= (0,0,1), e_4=( -1,2,2), e_5 = (-2,3,6),\cdots$.
  Note that $e_2, e_4, e_6, \cdots , e_{2n} $ take the top half of the graph and
   $e_1, e_3, e_5, \cdots , e_{2n+1}$  take the bottom half of the graph.

\medskip

\bigskip
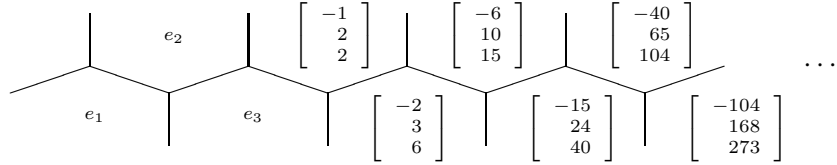
\begin{center}
 \begin{picture}(200,20)

\multiput(-30,0)(0,0){1}{\line(0,1){20}}
\multiput(-30,0)(0,0){1}{\line(-3,-1){30}}
\multiput(-30,0)(0,0){1}{\line(3,-1){30}}

\multiput(0,-10)(0,0){1}{\line(0,-1){20}}
\multiput(0,-10)(0,0){1}{\line(3,1){30}}

\multiput(30,0)(0,0){1}{\line(3,-1){30}}
\multiput(30,0)(0,0){1}{\line(0,1){20}}

\multiput(60,-10)(0,0){1}{\line(0,-1){20}}
\multiput(60,-10)(0,0){1}{\line(3,1){30}}


\multiput(90,0)(0,0){1}{\line(3,-1){30}}
\multiput(90,0)(0,0){1}{\line(0,1){20}}

\multiput(120,-10)(0,0){1}{\line(3,1){30}}
\multiput(120,-10)(0,0){1}{\line(0,-1){20}}

\multiput(150,0)(0,0){1}{\line(0,1){20}}
\multiput(150,0)(0,0){1}{\line(3,-1){30}}

\multiput(180,-10)(0,0){1}{\line(0,-1){20}}
\multiput(180,-10)(0,0){1}{\line(3,1){30}}

\put( -2, 10){\tiny$e_2$}
\put( -32, -20){\tiny$e_1$}

\put( 28, -20){\tiny$e_3$}

\put( 47, 10){\tiny$
\left [\begin{array} {r}
-1\\2\\2\\
\end{array}\right ]
$}

\put( 75, -25){\tiny$
\left [\begin{array} {r}
-2\\3\\6\\
\end{array}\right ]
$}

\put( 105, 10){\tiny$
\left [\begin{array} {r}
-6\\10\\15\\
\end{array}\right ]
$}

\put( 135, -25){\tiny$
\left [\begin{array} {r}
-15\\24\\40\\
\end{array}\right ]
$}

\put( 165, 10){\tiny$
\left [\begin{array} {r}
-40\\65\\104\\
\end{array}\right ]
$}

\put( 195, -25){\tiny$
\left [\begin{array} {r}
-104\\168\\273\\
\end{array}\right ]
$}

\put (240, -2) {$\cdots$}

\end{picture}
\end{center}

\vspace{1.5cm}

\noindent One sees immediately the following :

\begin{enumerate}
\item[(a)] The entries of the vectors $e_n =(a,b,c)$ are product of  two Fibonacci numbers.
To be more accurate, for the first nine terms,  the vectors take the form
$$(-F_nF_{n+1}, F_{n}F_{n+2}, F_{n+1}F_{n+2}).\eqno(5.3)$$
\item[(b)] The norm of $e_n= (a,b,c)$ is just the sum $a+b+c$,
  where the norm $N((x_1,x_2,x_3))$ is defined to be $(x_1^2+x_2^2+x_3^2)^{1/2}$
   (see Appendix A).
\item[(c)] The absolute value of the entries of $e_{2n}-e_{2n-2}$ (the top half of the
 trivalent graph) and
 $e_{2n+1}-e_{2n-1}$ (the bottom half of the trivalent graph) are Fibonacci numbers.
  For instance,
 ($-40, 65,104)-(-6,10,15) = (-F_9, F_{10}, F_{11})$. This allows one to write
  each entry of the vectors as a sum of Fibonacci numbers.

\end{enumerate}

\noindent
(a) and (b) of the above implies that for $n\le 6$
{\small $$(F_nF_{n+1})^2 + (F_nF_{n+2})^2 +(F_{n+1}F_{n+2})^2 = (-F_nF_{n+1} +F_nF_{n+2}+F_{n+1}F_{n+2})^2,
\eqno(5.4)$$}
\noindent  which leads us to (i) of the following lemma.
  Note that  $-F_nF_{n+1} +F_nF_{n+2}+F_{n+1}F_{n+2} = F_n^2 + F_{n+1}F_{n+2} $.
A careful study of (a) and (c) implies that each entry of the nine vectors can be written as
 sum as well as product of Fibonacci numbers which leads us to (ii)-(v) of the following lemma.

\medskip
\noindent
{\bf Lemma 5.1.} {\em Let $F_n$ be the $n$-th Fibonacci number. Then the following hold.}
\begin{enumerate}
\item[(i)]
$(F_nF_{n+1})^2 + (F_nF_{n+2})^2 +(F_{n+1}F_{n+2})^2 = (  F_n^2+ F_{n+1}F_{n+2}   )^2.$
\item[(ii)] $F_{2n-3}F_{2n-2} = F_1 +F_5 + \cdots + F_{4n-7}$.
\item[(iii)]$F_{2n-3}F_{2n-1}=  1+ F_2 +F_6 + \cdots +F_{4n-6}.$
\item[(iv)]$F_{2n-2}F_{2n-1} = F_3 +F_7 + \cdots +F_{4n-5}.$
\item[(v)] $F_{2n-2}F_{2n}=   F_4 +F_8 + \cdots+ F_{4n-4}.$
\end{enumerate}
\subsection{Discussion.}  (i) of the above lemma can be viewed as generalisation of Raine's results
 on Pythagorean's triple. (i)-(v) must be well known. As they are not included in
  [L] or  [W], we have them here for the reader's reference. Proofs of (i)-(v) are not
   included here as they can be proved easily.
    As the identities in Lemma 5.1 are from our observation of the trivalent graph.
     We consider that the recurrence relation (1.1), (5.1) and the
      trivalent graph $F(e_1,e_2,e_3)$ make those
      identities {\em visible}.

 \medskip
 \noindent To one's surprise, the trivalent graph actually tells us more.

    \begin{enumerate}
    \item[(i)] The sum of the first  entries (staring from $e_4$) of the first $2k-1$ consecutive
     vectors is the negative of a perfect square of a Fibonacci numbers.
 \item[(ii)] The sum of the second entry (starting from $e_4$) of the first $k$ vectors is
  a product of two Fibonacci numbers.
 \item[(iii)] The entries of every vector is a product of two
  Fibonacci numbers. Let
   $(-a,b,c)$ be such a vector.  Then $c-b-a = \pm 1$.

 \item[(iv)] Take {\em  any}  two consecutive vectors of the top half of the trivalent graph (such as
  $(e_2, e_4), $ $(e_4, e_6),\cdots )$. Label them as $(-a,b,c)$ and $(-A,B, C)$. Then
   $C-c = (B-b) + (A-a)$.
\item[(v)] Take two consecutive vectors of the top half (likewise
  the bottom half) and label them as $(-a,b,c)$ and $(-C,B,A)$. One {\em sees}  that
   all the entries are product
   of two Fibonacci numbers and the product of $a$ and $A$ is one less than a fourth power of a Fibonacci
    number ! ($1\cdot 15 = 2^4-1$, $6\cdot 104 = 5^4-1$, $2\cdot 40= 3^4 -1$, $15\cdot 273= 12^4-1$).
 \end{enumerate}

\noindent (i)-(v) of the above actually give us five well known identities. Take (v) for
     instance,  our observation shows that
     {\em
     a fourth power of a Fibonacci number $-1 = $ product of four Fibonacci numbers.}
To be more accurate, one has the
      remarkable Gelin-Ces\`{a}ro identity {\em visible}.
 $$F_n^4 -F_{n-2}F_{n-1}F_{n+1}F_{n+2} =1.\eqno(5.5)$$
 We are currently investigate  the trivalent graph
  $F(u, v,w)$ for arbitrary triples $(u,v.w)$. It turns out that such study  makes a lot of identities {\em visible}. For example, every identity appears in the right hand side of our table
   ((d1)-(d11))  in section 4 can be
   seen from some trivalent graph $F(u,v, w)$.
  See [LL]
   for more detail.
\section {Discussion}
We have demonstrated in this article that a simple
 study of the recurrence relation (1.1) ends up
  with a unified proof for  many  known identities
  in the literature. This suggests that one should probably group the
   identities together based on the recurrence relations (if it exists) and
    study them as a whole. Note that a given function may satisfy more than one recurrence
     relations ($(-1)^nF_n$ satisfies (i) of the following  and (3.5)).
   The next recurrence relation in line, we believe,
   should be
    \begin{enumerate}
    \item[(i)] $x(n+2) = -x(n+1)+ x(n)$,
    \item[(ii)] $x(n+3) = -2x(n+2) +2x(n+1) +x(n)$.
    \end{enumerate}
Identities (in Fibonacci numbers) with such  recurrence relations are rare but of great importance. To see our point,
 one recall that the right hand side of the very elegant identity of Melham's
 ($F_{n+1}F_{n+2}F_{n+6} -F_{n+3}^3=(-1)^{n}F_{n}$)
  satisfies
  (i)  of the above
  and that the following attractive identities of
  Fairgrieve and Gould ([FG])
  also satisfy (i) and (ii) of the above.
  $$F_{n-2}F_{n+1}^2 -F_n^3 = (-1)^n F_{n-1},\eqno(6.1)$$
   $$F_{n-3}F_{n+1}^3-F_n^4 = (-1)^n (F_{n-1}F_{n+3}+2F_n^2).\eqno(6.2)$$
To end our discussion, we give the following example which suggests how a new identity
 can be obtained by the study of recurrence relation (i):
     Since the right hand side of (6.1) satisfies (i) of the
     above, $x(n) = F_{n-2}F_{n+1}^2 -F_n^3$ satisfies the same recurrence relation.
      Namely, $x(n+2) = -x(n+1) +x(n)$. With the help of the famous identity
      $F_{3n} = F_{n+1}^3 + F_n^3 - F_{n-1}^3$, one has
      $$F_nF_{n+3}^2 +F_{n-1}F_{n+2}^2  -F_{n-2}F_{n+1}^2
      =F_{3n+3},\eqno(6.3)$$


\section{Appendix A}
 Let
 $u = (u_i) ,v=(v_i) ,w=(w_i) \in \Bbb Z^3 $ be vectors.
 We say $\{u,v, w\}$ is a $\mathcal F$-triple if
 $N(u)=
 u_1+u_2+u_3, N(v)
= v_1+v_2+v_3
$ and  $N(w)=w_1+w_2+w_3$ are squares in $\Bbb N$ and
\begin{enumerate}
\item[(i)]
$2u\cdot v - v\cdot w - w\cdot u
 = 2N(u)N(v)  - N(v)N(w) - N(w)N(u)$,
\item[(ii)]
$2u\cdot w - v\cdot w - v\cdot u
 =   2N(u)N(w)  - N(v)N(w) - N(v)N(u)  $,
\item[(iii)]
$2v\cdot w - v\cdot u - w\cdot u
 = 2N(v)N(w)  - N(v)N(u) - N(w)N(u)$,
\end{enumerate}
where $u\cdot v$ is the usual dot  product.
One sees easily that $\{e_1, e_2, e_3\}$ is an $\mathcal F$-triple,
 where $e_1=(1,0,0), e_2= (0,1,0), e_3=(0,0,1)$. The following lemma shows that
  if $\{u, v, w\}$ is an $\mathcal F$-triple, then any vector $(a,b,c)$ in $F(u,v, w)$
   has the property $N((a,b,c)) = a+b+c$, which proves (b) of section 5.

\medskip
\noindent {\bf Lemma A.} {\em Let
 $u = (u_i) ,v=(v_i) ,w=(w_i) \in \Bbb R^3 $ be a $\mathcal F$-triple
and let  $x = (x_i) =2(u+v)-w$,
$y =(y_i) =2(w+v)-u$,
$z =(z_i) =2(w+u)-v$.
Then  the following hold.}
\begin{enumerate}
\item[(i)]
$\{u,v, x \}$, $ \{u,w, z\}$ and
 $\{v,w, y\}$ are $\mathcal F$-triples,
\item[(ii)]
$N(x)^2 = (2N(u) + 2N(v)-N(w))^2$,
$N(y)^2 = (2N(w) + 2N(v)-N(u))^2$,
$N(z)^2 = (2N(w) + 2N(u)-N(v))^2$.
\end{enumerate}

\smallskip
\noindent
 {\em Proof.}  The lemma is straightforward and can be proved by direct calculation.
 \qed

 \medskip
\noindent Note that  $x$, $y$, and $z$ in the above lemma are defined as in $(5.1)$ and can be
 described as follows :
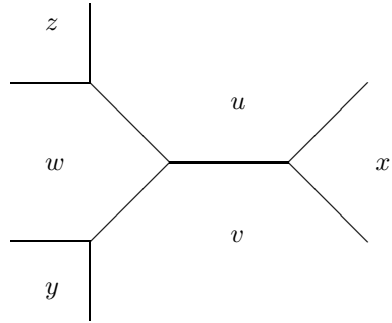
\begin{center}
 \begin{picture}(60,75)

 \multiput(2,0)(30,0){1}{\line(1,0){45}}

\multiput(2,0)(60,0){1}{\line(-1,1){30}}
 \multiput(2,0)(60,0){1}{\line(-1,-1){30}}

 \multiput(77,-30)(60,0){1}{\line(-1,1){30}}
 \multiput(77,30)(30,30){1}{\line(-1,-1){30}}

 \multiput(-28,-30)(50,0){1}{\line(-4,0) {30}}
 \multiput(-28,-60)(50,0){1}{\line(0,4){30}}

  \multiput(-58,30)(50,0){1}{\line(4,0){30}}
 \multiput (-28,30)(50,0){1}{\line(0,4){30}}

\put(80 , -3){\small  $x$}
\put(25, 20){\small $u$}
\put(25, -30){\small $v$}

\put(-45 , -3){\small  $w$}
\put(-45 , 50){\small  $z$}
\put(-45 , -50){\small  $y$}

 \end{picture}
\end{center}

\vspace{2.5cm}

\bigskip
\noindent Following our lemma, one may extend the above graph to an infinite trivalent
graph that takes the whole $xy$-plane such that
each triple $\{r,s,t\}$ associated to a vertex is an $\mathcal F$-triple. In particular,
the entries of every vector of this trivalent graph give solution to
  $x^2+y^2 +z^2 = (x+y+z)^2$. Note that a complete set of integral solutions of the above mentioned equation
   is given by $\{ (mn, m(m+n), n(m+n)) \,:\, n, m \in \Bbb Z\}$.

\section {Appendix B : More Recurrence relations}
Let $x(n)$ be a function defined on $\Bbb Z$. Consider the equation
$$ x(n+k) = a_{k-1} x(n+k-1) +\cdots + a_1x(n+1) + a_0x(n).\eqno(B1)$$
 One sees easily that whether $x(n)$ satisfies some recurrence relation depends on whether
  there exists some $k$ and $a_i$'s such that $(B1)$ holds for all $n$. In the case $x(n)$ indeed
   admits some recurrence relation, such relation can be obtained by solving system of
    linear equations.

\subsection{The Recurrence relation
 $x(n+4) = 3x(n+3)+6x(n+2)-3x(n+1) -x(n)$}
In [M2], Melham proved that
$$F_{n+1}F_{n+2}F_{n+6} -F_{n+3}^3=(-1)^{n}F_{n}.\eqno(B2)$$

\medskip
\noindent We shall give our alternative proof as follows.
Let $A(n) = F_{n+1}F_{n+2}F_{n+6} -F_{n+3}^3$, $B(n) = (-1)^{n}F_{n}$.
One sees easily that both $A(n)$ and $B(n)$ satisfy the above recurrence relation.
 Since
 \begin{enumerate}
 \item[(i)]  $A(n)$ and $B(n)$ satisfy the above recurrence relation and $A(n) = B(n) $ for $n = 0,1,2,3$,
 \item[(ii)] function $x(n)$ satisfies the above recurrence relation  is completely
  determined by $x(0)$,  $x(1),$
  $ x(2)$ and $x(3)$,
  \end{enumerate}
  we conclude that $A(n) = B(n)$. This completes the proof of $(B2)$.
The identity $F_{3n} = F_{n+1}^3 +F_n^3-F_{n-1}^3$ and
 Fairgrieve and Gould's identities
  ((11), (12)  of [FG]) can be proved by the same method.

\subsection{Recurrence relation for $F_n^4$} $F_n^4$ satisfies the following recurrence
 relation,
 $$x(n+5)
  = 5x(n+4)+15x(n+3)-15x(n+2)-5x(n+1)+x(n).\eqno(B3)$$
One sees easily that both the left and right hand side of $(6.2)$ satisfy $(B3)$.
 As a consequence,  identity (6.2) can be verified by applying our technique given in
  the above subsection.
\subsection{Construction of Identities} Recurrence relations can be used to construct 
 identities. Take (6.3) for example, one can actually construct (6.3) as follows. 
 
 $$\begin{array}{lrrrr}
 
 x(n) & x(0) &x(1)& x(2)& x(3)\\
 \\
 F_{3n+3} & 2 &8 & 34& 144\\
 
 \\
  F_nF_{n+3}^2 & 0  & 9 &25  &128 \\
  
  \\
  F_{n-1}F_{n+2}^2& 1 & 0& 9 &  25\\
  \\
  F_{n-2}F_{n+1}^2 & -1 & 1 & 0 & 9 \\
\\
\end{array}
$$
Since $F_{3n+3},F_nF_{n+3}^2,F_{n-1}F_{n+2}^2$ and $  F_{n-2}F_{n+1}^2 $ satisfy 
 the recurrence relation
 $x(n+4) = 3x(n+3)+6x(n+2)-3x(n+1) -x(n)$, applying (ii) of the 
  above, one sees from the above table that 
  $$F_{3n+3}=F_nF_{n+3}^2+F_{n-1}F_{n+2}^2 -   F_{n-2}F_{n+1}^2.\eqno(B4)$$

\medskip

\noindent MSC2010: 11B39, 11B83


\bigskip






\begin{thebibliography}{99}

\bibitem[FG]{sloane}
S. Fairgrieve, H. W. Gould \emph{Product difference Fibonacci identities of Simson,
 Gelin-Ces\`{a}ro, Tagiuri and generalizations}, The Fibonacci
Quarterly \textbf{43.2} (2005), 137-141.


\bibitem[HB]{sloane}
V. E. Hoggatt, Jr, M. Bicknell, \emph{Some new Fibonacci identities}, The Fibonacci
Quarterly \textbf{2.1} (1964), 29--32.

\bibitem[Hor]{sloane}
A. F. Horadam, \emph{Basic Properties of a Certain Generalized Sequence of Numbers}, The Fibonacci
Quarterly \textbf{3.3} (1965), 161--76.





\bibitem[How]{sloane}
F. T. Howard, \emph{The sum of the squares of two generalized Fibonacci numbers}, The Fibonacci
Quarterly \textbf{41.1} (2003), 80--84.



\bibitem[L]{sloane}
C. T. Long , \emph{Discover Fibonacci identities}, The Fibonacci
Quarterly \textbf{24.2} (1986), 160--167.



\bibitem[LL]{sloane}
C.L. Lang and M.L. Lang, \emph{Fibonacci Numbers and Trivalent graphs}, preprint.



\bibitem[M1]{sloane}
R.S Melham, \emph{Families of identities involving sums of powers of the Fibonacci and Lucas numbers}, The Fibonacci
Quarterly \textbf{37.4} (1999), 315--319.


\bibitem[M2]{sloane}
R.S Melham, \emph{A Fibonacci identity in the spirit of Simon and
 Gelin-Ces\`{a}ro}, The Fibonacci
Quarterly \textbf{41.2} (2003), 142-143.






\bibitem[W]{sloane} E. W.  Weisstein,
{\em Fibinacci numbers,}
{\tt http:// www.mathworld.wolfram.com/FibonacciNumbers.html}.



\end{thebibliography}
\end{document}